\documentclass[12pt]{amsart}

\usepackage{amsmath,amsthm}
\usepackage{amssymb,latexsym}
\usepackage{graphicx}
\theoremstyle{definition}

\theoremstyle{remark}

\numberwithin{equation}{section}

\newcommand{\To}{\longrightarrow}

\newcommand{\bee}{\begin{equation}}
\newcommand{\ene}{\end{equation}}
\newcommand{\pa}{\partial}
\newcommand{\la}{\lambda}
\newcommand{\al}{\alpha}

\begin{document}

\title[The utilization of total mass]{The utilization of total mass to determine the switching points in the symmetric boundary control of a diffusion problem}%

\author{M. Salman}
\address{Department of Math, Physics \& Stat,
 Qatar University, PO Box 2713, Doha, Qatar}
\email{msalmanz@gmail.com}


\begin{abstract}
The authors study the problem $u_t=u_{xx},\ 0<x<1,\ t>0; \
u(x,0)=0,$ and $u(0,t)=u(1,t)=\psi(t),$ where $\psi(t)=u_0$ for
$t_{2k} < t<t_{2k+1}$ and $\psi(t)=0$ for $t_{2k+1} <t<t_{2k+2},\
k=0,1,2,\ldots$ with $t_0=0$ and the sequence $t_{k}$ is
determined by the equations $\int_0^1 u(x,t_k)dx = M,$ for
$k=1,3,5,\dots,$ and $\int_0^1 u(x,t_k)dx = m,$ for
$k=2,4,6,\dots$ and where $0<m<M<u_0$. Note that the switching
points $t_k,\quad k=1,2,3,\ldots$ are unknown. Existence and
uniqueness are demonstrated. Theoretical estimates of the $t_k$
and $t_{k+1}-t_k$ are obtained and numerical verifications of the
estimates are presented.

\end{abstract}
\maketitle


\numberwithin{equation}{section}








\section{Introduction}

 As motivation for the mathematical problems considered
in this work, consider a chamber in the form of a long linear
transparent tube. We allow for the introduction or removal of
material in a gaseous state at the ends of the tube. The material
diffuses throughout the tube with or without reaction with other
materials. By illuminating the tube on one side with a light source
with a frequency range spanning the absorption range for the
material and collecting the residual light that passes through the
tube with photo-reception equipment, we can obtain a measurement of
the total mass of material contained in the tube as a function of
time. Using the total mass as switch points for changing the
boundary conditions for introduction or removal of material. The
objective is to keep the total mass of material in the tube
oscillating between two set values such as $m<M$. The physical
application for such a system is the control of reaction diffusion
systems such as production of a chemical material in a reaction
chamber via the introduction of reactants at the boundary of
chamber.

In this work we study the diffusion equation \bee\label{1.1}
u_t=u_{xx}, \quad 0<x<1, \;\; t>0 \ene subject to the initial
concentration
$$
u(x,0)=0, \quad 0<x<1
$$
and boundary conditions controlled by the total mass $\mu
(t)=\int_0^1 u(x,t)dx$.  We are going to begin by setting the
concentration to be $u=u_0$ at both boundary points $x=0$ and
$x=1$, where $u_0$ is a positive constant. We shall watch the
total mass $\int_0^1 u(x,t)dx$ until it reaches a certain
specified level $M$ at time $t=t_1$, where $0<M<u_0$. At this
moment, $t=t_1$ we switch the concentration to
$$
u(0,t)=u(1,t)=0 \quad t_1<t.
$$
We keep watching the total mass until it drops down to a
prespecified level $m$ at time $t=t_2$, where $0<m<M<u_0$. We keep
switching the concentration according to the level of the total
mass so that we always have
$$
m\leq \int_0^1 u(x,t)dx\leq M.
$$
In other words, the boundary conditions will be \bee\label{1.2}
u(0,t)=u(1,t)=:\phi (t)= \begin{cases}
u_0, & t_{2n}\leq t\leq t_{2n+1},\\
0, & t_{2n+1} \leq t\leq t_{2n+2} \end{cases} \ene where
$n=0,1,2,\dots$; and the sequence $\{t_n\}$ will be strictly
increasing, i.e.
$$
0=t_0<t_1<t_2<\dots
$$
and its terms are defined by the equations
\begin{align*}
&\int_0^1 u(x,t_n)dx = M, \quad n=1,3,5,\dots \; ,\\
&\int_0^1u(x,t_n)dx = m, \quad n=2,4,6,\dots \; .
\end{align*}

\section{Existence of the Sequence $\{t_n\}$}

For the sake of simplicity, we will take $u_0=1$. Consider the
problem
\begin{align}\label{2.1}
&u_t =u_{xx}, \quad 0<x<1, \quad t>0,\notag\\
&u(0,t) =u(1,t)=\phi (t), \quad t>0,\notag\\
&u(x,0) = 0 .
\end{align}
Then the solution will be (see \cite{can1})  \bee\label{2.2}
u(x,t)= \int_0^t \left[ \frac{\pa\theta}{\pa x} (x-1,t-\tau
)-\frac{\pa\theta}{\pa x} (x,t-\tau )\right] \phi (\tau )d\tau,
\ene where \cite{kev} \bee\label{2.3} \theta (x,t)=\frac12
+\sum_{n=1}^\infty e^{-n^2\pi^2t} \cos n\pi x. \ene Upon
integrating (\ref{2.2}) with respect to $x$ from 0 to 1 we get
\bee\label{2.4} \mu(t):=\int_0^1 u(x,t)dx=4\int_0^t [\theta
(0,t-\tau )-\theta (1,t-\tau )]\phi (\tau )d\tau . \ene
Substituting (\ref{2.3}) in (\ref{2.4}), we obtain \bee\label{2.5}
\mu(t)=8\int_0^t \phi(\tau ) \left(\sum_{k=0}^\infty
e^{-\la_{2k+1}^2 (t-\tau )} \right) d\tau \ene where $\la_k=k\pi$.

The first stage we set $\phi (t)=1$. This will give
$$
\mu(t)=8\sum_{k=0}^\infty \frac{1}{\la_{2k+1}^2} \left[
1-e^{-\la_{2k+1}^2 t} \right]
$$
which is a strictly increasing function of $t$, and it ranges
between 0 and 1 as $t$ ranges from 0 to $\infty$. Therefore, there
exists a positive $t_1$ such that
$$
\mu(t_1)=M, \quad 0<M<1.
$$
For $t>t_1$, we set $\phi (t)=0$. Then equation (\ref{2.5})
implies
\begin{align*}
\mu(t)&=8\int_0^{t_1} \sum_{k=0}^\infty e^{-\la_{2k+1}^2(t-\tau )} d\tau , \quad t>t_1\\
&= 8\sum_{k=0}^\infty \frac{e^{-\la_{2k+1}^2t}}{\la_{2k+1}^2}
\left[ e^{\la_{2k+1}^2t_1} -1\right] ,
\end{align*}
which is a strictly decreasing function that falls from $M$ to 0
as $t$ goes from $t_1$ to $\infty$. Hence, there exists a $t_2$
such that
$$
\mu (t_2)=m,\quad0<m<M<1.
$$
In an inductive fashion, we obtain as $t>t_{2n}$ and $\phi (t)=1$,
\begin{align*}
\mu(t) &=8\sum_{j=0}^{n-1} \int_{t_{2j}}^{t_{2j+1}} \sum_{k=0}^\infty e^{-\la_{2k+1}^2(t-\tau )} d\tau \\
&\quad +8\int_{t_{2n}}^t \sum_{k=0}^\infty e^{-\la_{2k+1}^2(t-\tau )}d\tau , \quad t>t_{2n}\\
&= 8\sum_{k=0}^\infty \frac{e^{-\la_{2k+1}^2t}}{\la_{2k+1}^2} \sum_{j=0}^{n-1} \left[ e^{\la_{2k+1}^2t_{2j+1}}-e^{\la_{2k+1}^2t_{2j}} \right]\\
&\quad +8\sum_{k=0}^\infty \frac{1}{\la_{2k+1}^2} \left[ 1-e^{-\la_{2k+1}^2(t-t_{2n})}\right] \\
&=8\sum_{k=0}^\infty \frac{1}{\la_{2k+1}^2} \left\{
1-e^{-\la_{2k+1}^2t} \left[ \sum_{j=0}^{2n} (-1)^j
e^{\la_{2k+1}^2t_j} \right] \right\}
\end{align*}
which increases from $m$ to 1 as $t$ goes from $t_{2n}$ to
infinity. Thus, there exists a $t_{2n+1}$ such that
$$
\mu(t_{2n+1})=M.
$$
When $t>t_{2n+1}, \phi(t)=0,$ which implies
\begin{align*}
\mu(t) &=8\int_0^{t_{2n+1}}\phi(\tau )\sum_{k=0}^\infty e^{-\la_{2k+1}^2(t-\tau )} d\tau , \quad t>t_{2n+1}\\
&= 8 \sum_{j=0}^n \int_{2_{2j}}^{2j+1} \sum_{k=0}^\infty e^{-\la_{2k+1}^2(t-\tau )} d\tau\\
&=8\sum_{k=0}^\infty \frac{e^{-\la_{2k+1}^2t}}{\la_{2k+1}^2}
\left\{ \sum_{j=0}^{2n+1} (-1)^{j+1} e^{\la_{2k+1}^2 t_j} \right\}
.
\end{align*}
Hence, $\mu(t)$ will continuously decrease down from $M$ to 0 as
$t$ goes from $t_{2n+1}$ to infinity. This ensures the existence
of $t_{2n+2}$ such that
$$
\mu(t_{2n+2})=m.
$$
From the argument above, we have constructed the sequence
$\{t_n\}$, where
$$
0=t_0<t_1<t_2<\dots \; .
$$
Given the switching sequence $\{t_n\}$, the existence and
uniqueness of the solution follows immediately.

\section{A First Term Approximation}

In this section we study eqnuation (\ref{2.5}) by using the first
term of the infinite series, that is \bee\label{3.1} \mu(t) \simeq
8\int_0^t \phi(\tau )e^{-\la_1^2(t-\tau )} d\tau =: \tilde \mu (t)
\ene where $\la_1=\pi$.

We will find an approximation to the sequence $\{t_n\}$ in an
explicit form. For the sake of simplicity, the upper and lower
bound on total mass are taken to be \bee\label{3.2}
M=\frac{8}{\pi^2}(1-\al ) \ene \bee\label{3.3} m=\frac{8}{\pi^2}
\al \quad\quad\; \ene where the restriction $ 0<\al <\frac12$ has
been set in order to $0<m<M$.  The first time step $t_1$ can be
calculated through the equation
$$
\tilde \mu (t_1)=M,
$$
which gives,
$$
\frac{8}{\la_1^2} \left[ 1-e^{-\la_1^2t_1}\right]
=\frac{8}{\la_1^2} (1-\al ),
$$
i.e. \bee\label{3.4} e^{-\la_1^2t_1}=\al. \ene  The second time
step can be found by
$$
\tilde \mu (t_2)=m.
$$
This gives

\bee \frac{8e^{-\la_1^2t_2}}{\la_1^2} \left[
e^{\la_1^2t_1}-1\right] =\frac{8}{\pi^2} \al .\ene By using
(\ref{3.4}), the above equation implies \bee\label{3.6}
e^{-\la_1^2t_2} =\frac{\al^2}{1-\al} .\ene Using a similar
argument, we can inductively obtain

$$ e^{-\la_1^2t_3} =\frac{\al^3}{(1-\al)^2} ,$$
$$e^{-\la_1^2t_4} =\frac{\al^4}{(1-\al)^3}, \\$$
$$\vdots\\$$

\bee\label{3.8} e^{-\la_1^2t_n} =\frac{\al^n}{(1-\al)^{n-1}},\quad
n\ge1. \ene

Next we use mathematical induction to prove (\ref{3.8}). The point
$t_{2n}$ can be obtained by
\begin{align*}
\tilde \mu(t_{2n}) &= 8\int_0^{t_{2n-1}} \phi(\tau )e^{-\la_1^2(t_{2n}-\tau )}d\tau\\
&= \left[ \int_0^{t_1} e^{-\la_1^2(t_{2n}-\tau )} d\tau  +\dots +\int_{t_{2n-2}}^{t_{2n-1}} e^{-\la_1^2(t_{2n}-\tau )}d\tau \right]\\
&=\frac{8e^{-\la_1^2t_{2n}}}{\la_1^2} \left[ e^{\la_1t_{2n-1}}
-e^{\la_1t_{2n-2}} +\dots +e^{\la_1t_1} -1\right] =m .
\end{align*}
Assuming formula (\ref{3.8}) is valid for all $k<2n$, then the
last equation is equivalent to
$$
e^{-\la_1^2t_{2n}} \left[\frac{(1-\al)^{2n-2}}{\al^{2n-1}}
-\frac{(1-\al )^{2n-3}}{\al^{2n-2}} +\dots +\frac{1}{\al}
-1\right] =\al.
$$
Upon adding all the terms inside the brackets as a finite
geometric sum, we will obtain \bee\label{3.8e}
e^{-\la_1^2t_{2n}}=\frac{\al^{2n}}{(1-\al )^{2n-1}}. \ene The time
step $t_{2n+1}$ can be found as a solution of
\begin{align*}
\tilde \mu (t_{2n+1}) &=8\int_0^{t_{2n+1}} \phi(\tau )e^{-\la_1^2(t_{2n+1}-\tau )}d\tau\\
&= 8\left[ \int_0^{t_1}+\int_{t_2}^{t_3} +\dots +\int_{t_{2n}}^{t_{2n+1}} e^{-\la_1^2(t_{2n+1}-\tau )}d \tau \right]\\
&=\frac{8e^{-\la_1^2t_{2n+1}}}{\la_1^2} \left[e^{\la_1^2t_{2n+1}}-e^{\la_1^2t_{2n}} +e^{\la_1^2t_{2n-1}}  -e^{-\la_1 t_{2n-2}} + \dots + e^{\la_1^2t_1} -1 \right] \\
&=M.
\end{align*}
Assuming the validity of (\ref{3.8}) and using (\ref{3.8e}) for
all $k\leq 2n$, the above equation implies
$$
e^{-\la_1^2t_{2n+1}} \left[ \frac{(1-\al )^{2n-1}}{\al^{2n}}
-\frac{(1-\al )^{2n-2}}{\al^{2n-1}} +\dots -\frac{1}{\al}+1\right]
=\al .
$$
If we add these terms which constitute finite geometric series, we
can obtain \bee\label{3.9} e^{-\la_1^2t_{2n+1}}
=\frac{\al^{2n+1}}{(\al-1)^{2n}} . \ene Hence, formula (\ref{3.8})
is valid.\\ Formula (\ref{3.8}) implies that
$$
t_n=\frac{1}{\pi^2} \ln \frac{(1-\al )^{n-1}}{\al^n}, \quad n\geq
1.
$$
Therefore,
$$
t_{n+1}-t_n=\frac{1}{\pi^2} \ln\frac{1-\al}{\al}, \quad n\geq 1,
$$
which is obviously independent of $n$.

\section{A Higher Order Approximation}

In this section we will study a higher order approximation of the
sequence $\{ t_n\}$. This approximation is based on dropping all
the terms that has $e^{-\lambda_{2k+1}^2t_n}$ from the expression
of the mass, except for a few dominating terms. The upper and
lower bounds on the mass are not necessarily symmetric. For
simplicity of computation, the bounds are chosen to be
$$
M=8\left[ \sum_{k=0}^\infty \frac{1}{\lambda_{2k+1}^2}
-\frac{\al}{\la_1^2}\right]
$$
and
$$
m=\frac{8}{\la_1^2} \beta
$$
where $\la_{2k+1}^2=(2k+1)^2\pi^2, k=0,1,\dots$; $\al$ and $\beta$
are positive numbers that are chosen so that the inequality
$0<m<M$ holds. Since $\sum_{k=0}^\infty
\frac{1}{\la_{2k+1}^2}=\frac{1}{\pi^2}\sum_{k=0}^\infty
\frac{1}{(2k+1)^2} =\frac18$, then we have
$M=1-\frac{8}{\pi^2}\al$. \\To find $t_1$, we need to solve the
equation
$$
\mu(t_1)=M
$$
which is
$$
8\left[ \frac{1-e^{-\la_1^2t_1}}{\la_1^2} +\sum_{k=1}^\infty
\frac{1-e^{-\la_{2k+1}^2t_1}}{\la_{2k+1}^2} \right]=M.
$$
Upon dropping all the terms that have $e^{-\la_{2k+1}^2t_1}$ for
all $k\geq 1$, we get
$$
8\left[\sum_{k=0}^\infty \frac{1}{\la_{2k+1}^2}
-\frac{e^{-\la_1^2t_1}}{\la_1^2}\right] =1-\frac{8}{\pi^2}\al
$$
which implies \bee\label{4.1} e^{-\la_1^2t_1}=\al . \ene In a
similar fashion, we compute $t_2$ by solving
$$
\mu(t_2)=8\int_0^{t_1}\sum_{k=0}^\infty
e^{-\la_{2k+1}^2(t_2-\tau)}d\tau =m
$$
i.e.,
$$
8\sum_{k=0}^\infty \frac{e^{-\la_{2k+1}^2t_2}}{\la_{2k+1}^2}
\left[ e^{\la_{2k+1}^2t_1} -1\right] =m .
$$
By dropping all the terms that have $e^{\la_{2k+1}^2t_2}$ for all
$k\geq 1$, we get
$$
e^{-\la_1^2t_2} \left( e^{\la_1^2t_1} -1\right) =\beta .
$$
Using (4.1) we obtain \bee e^{-\la_1^2t_2} =\frac{\al\beta}{1-\al}
. \ene By employing a similar method of approximations we can find
$t_3$ through the equation
$$
e^{-\la_1^2t_3} \left[e^{\la_1^2t_2} -e^{\la_1^2t_1} +1\right]
=\al
$$
which implies that \bee e^{-\la_1^2t_3} =\frac{\al^2\beta}{(1-\al
)(1-\beta )} . \ene For $t_4$ and $t_5$, we obtain the explicit
expressions \bee e^{-\la_1^2t_4} =\frac{\al^2\beta^2}{(1-\al
)^2(1-\beta )} \ene and \bee e^{-\la_1^2t_5}
=\frac{\al^3\beta^2}{(1-\al )^2(1-\beta)^2} . \ene Thus, we choose
\bee e^{-\la_1^2t_{2n}} =\frac{\al^n\beta^n}{(1-\al )^n(1-\beta
)^{n-1}} , \quad n\geq 1 \ene and \bee e^{-\la_1^2t_{2n+1}}
=\frac{\al^{n+1}\beta^n}{(1-\al )^n(1-\beta )^n} , \quad n\geq 1 .
\ene By an induction argument similar to the one used in Section 3
we can show that (4.6) and (4.7) hold. Therefore the consecutive
time steps can be calculated by \bee t_{2n}-t_{2n-1} =\ln
\frac{1-\al }{\beta} \ene and \bee t_{2n+1}-t_{2n} =\ln
\frac{1-\beta }{\al} \ene where $n\geq 1$. For the special case
when $\al =\beta$, the time intervals reduce that shown at the end
of section 3.

\section{Numerical Results}

In this section, we use a finite difference technique along with
the trapezoidal rule to get an approximate discrete solution
$U_j^n$ and the sequence $\{ T_m\}$ when the total mass hits one
of the limits $M$ or $m$.

We discretize the space and time by using

i) $\Delta x=\frac{1}{J} , \quad X_j=j\Delta x, \quad j=0,1,\dots
,J$

ii) $\Delta t=\frac{T}{N} , \quad \tau_n =n\Delta t, \quad
n=0,\dots ,N$

\noindent where $J$ and $N$ are positive integer and $T$ is a
positive real number. The integer $N$ has to be chosen large
enough so that the time step $\Delta t$ is much smaller than the
differences $T_n-T_{n-1}$.

We consider the backward implicit finite difference scheme
$$
\frac{U_j^{n+1}-U_j^{n}}{\Delta t} = a \frac{U_{j-1}^{n+1}
-2U_j^{n+1} +U_{j+1}^{n+1}}{(\Delta x)^2}
$$
as a discretized version of $u_t=a u_{xx}$, where $a$ is a positive
constant. The above scheme can be written in the form \bee -b
U_{j-1}^{n+1} +(1+2b)U_j^{n+1} -b U_{j+1}^{n+1} =U_j^n \ene where
$j=1,\dots ,J-1$ and $n=0,1,\dots, N-1$. The initial data are set to
be $U_j^0=0$ for $j=1,\dots ,J-1$, and the boundary conditions are
$U_0^n =U_J^n =\phi(\tau_n )$ for $n=0,1,\dots ,N$. The function
$\phi (\tau_n )$ will be either 10 or 0 depending on the value of
the mass which will be approximated by the trapezoidal rule \bee
\mu_n =\frac{h}{2} \sum_{j=0}^{N-1} \left( U_j^{n+1} +U_{j+1}^n
\right). \ene The numerical experiment is carried out in the
following way. We start by setting the boundary conditions
$U_0^n=U_J^n =10$ then we solve a tridiagonal system coming out of
the difference method. We check the total mass $\mu_n$ in (5.2). We
keep doing that at each time step until the mass $\mu_n$ exceeds or
equals the upper limit $M$. Then we switch the boundary conditions
to $U_0^n=U_J^n=0$, and we continue the finite difference scheme for
several time steps $\Delta t$ until the total mass $\mu_n$ decreases
to $m$. At this moment we switch the boundary condition back to 10
and continue the process as we did before.

For the data specifications $\Delta x= 0.02$, $T=20$ $\Delta t=0.1$,
$a =0.05$, $M=7$, $m=3$, Table (1) shows the times switches $T_n$.
As we can see there, the duration of each stage
turns out to be constant.\\
For the same set of data, Graphs (1) through (6) show the
concentration versus the space. The graphs are obtained for
different stages, where at each stage the concentration is kept
constant at the end points.\\
 A profile of the concentrations at $x=0.5$ for various times is shown
  in Graph (7) with the same specified data.\\
For a different set of upper and lower bounds on the mass $M=5$
and $m=2$ along with $\Delta x=0.02,$ $T=20$, $\Delta t=0.02$,
Table (2) shows the time switches $T_n$. Durations of the
time intervals fluctuates between 0.5, 1.\\
Graph (8) is the concentration at $x=0.5$ for the same data
generating Table (2).

\noindent\textsc{ \underline {Conclusion:} } From Table 1 and
Table 2, we see that the theoretical estimates of
$t_{2n}-t_{2n-1}$ and $t_{2n+1}-t_{2n}$ are exihibited in the
numerical examples.

\begin{table}

\begin{tabular}{ | r | r | r |r | r | r |}
\hline
  n & Tn &  Tn - Tn-1 & n & Tn &  Tn - Tn-1\\\hline
  1 &  2.1000 & 2.1000 &  8 & 10.5000 & 1.2000\\
  2 &  3.3000 & 1.2000 &  9 & 11.7000 & 1.2000\\
  3 &  4.5000 & 1.2000 & 10 & 12.9000 & 1.2000\\
  4 &  5.7000 & 1.2000 & 11 & 14.1000 & 1.2000\\
  5 &  6.9000 & 1.2000 & 12 & 15.3000 & 1.2000\\
  6 &  8.1000 & 1.2000 & 13 & 16.5000 & 1.2000\\
  7 &  9.3000 & 1.2000 & 14 & 17.7000 & 1.2000\\\hline

\end{tabular}
\caption{ Time switches $T_n$ corresponding to $\Delta x=0.02$,
$T=20$, $\Delta t=0.1$, $a =0.05$, $M=7$, $m=3$ }

\end{table}

\begin{table}

\begin{tabular}{ | r | r | r |r | r | r |}
\hline
   n &  Tn       & Tn - Tn-1& n &  Tn       & Tn - Tn-1\\\hline

   1 &  1.0000  & 1.0000& 14 & 10.9400  & 1.0200\\
   2 &  1.9400  & 0.9400& 15 & 11.4200  & 0.4800\\
   3 &  2.4200  & 0.4800& 16 & 12.4400  & 1.0200\\
   4 &  3.4400  & 1.0200& 17 & 12.9200  & 0.4800\\
   5 &  3.9200  & 0.4800& 18 & 13.9400  & 1.0200\\
   6 &  4.9400  & 1.0200& 19 & 14.4200  & 0.4800\\
   7 &  5.4200  & 0.4800& 20 & 15.4400  & 1.0200\\
   8 & 6.4400   &1.0200 & 21 & 15.9200  & 0.4800\\
   9 &  6.9200  & 0.4800& 22 & 16.9400  & 1.0200\\
  10 &  7.9400  & 1.0200& 23 & 17.4200  & 0.4800\\
  11 &  8.4200  & 0.4800& 24 & 18.4400  & 1.0200\\
  12 &  9.4400  & 1.0200& 25 & 18.9200  & 0.4800\\
  13 &  9.9200  & 0.4800& 26 & 19.9400  & 1.0200\\\hline

\end{tabular}

\caption{}
\end{table}


\begin{figure}[b]
\scalebox{0.6}{\includegraphics{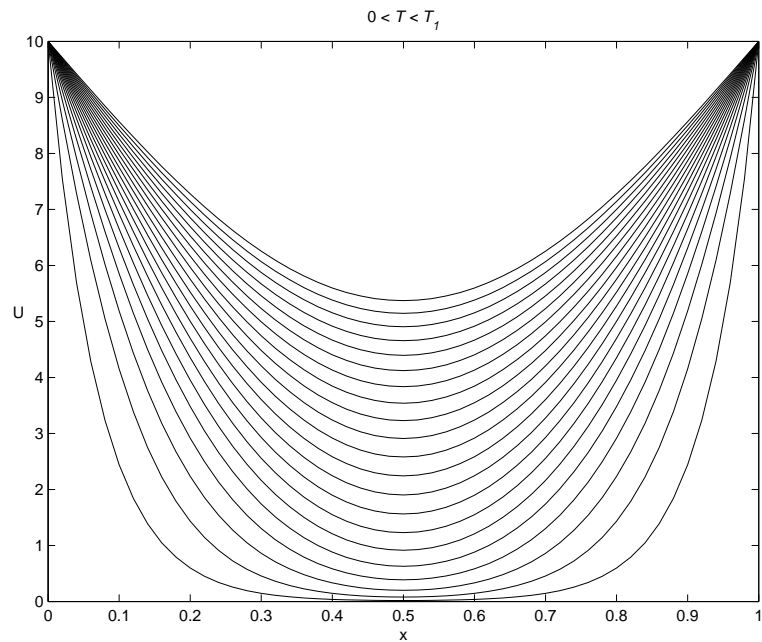}}
 \caption{The first
stage where the concentration $U$ is held at 10 at the end points.
Each curve shows the concentration profile at various discrete
time steps $t_n=n \Delta t$. As the time goes on, the level of
concentrations gets higher }
\end{figure}


\begin{figure}[b]
\scalebox{0.6}{\includegraphics{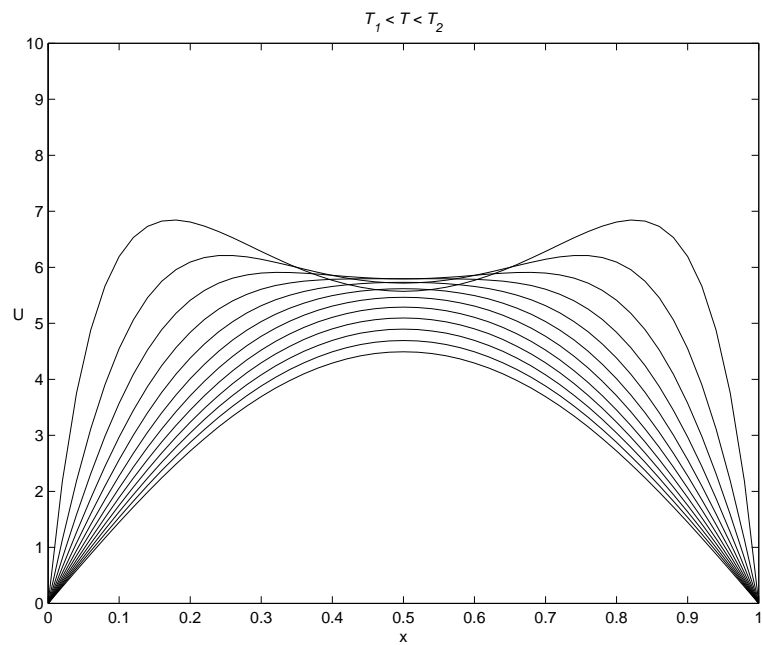}} \caption{The second
stage where the concentration $U$ is held at 0 at the end points.
 As the time goes on, the level of
concentrations, roughly speaking, decreases. Notice the
fluctuations when the concentration is dropped suddenly to 0 at
the beginning of the stage}
\end{figure}

\begin{figure}[b]
\scalebox{0.6}{\includegraphics{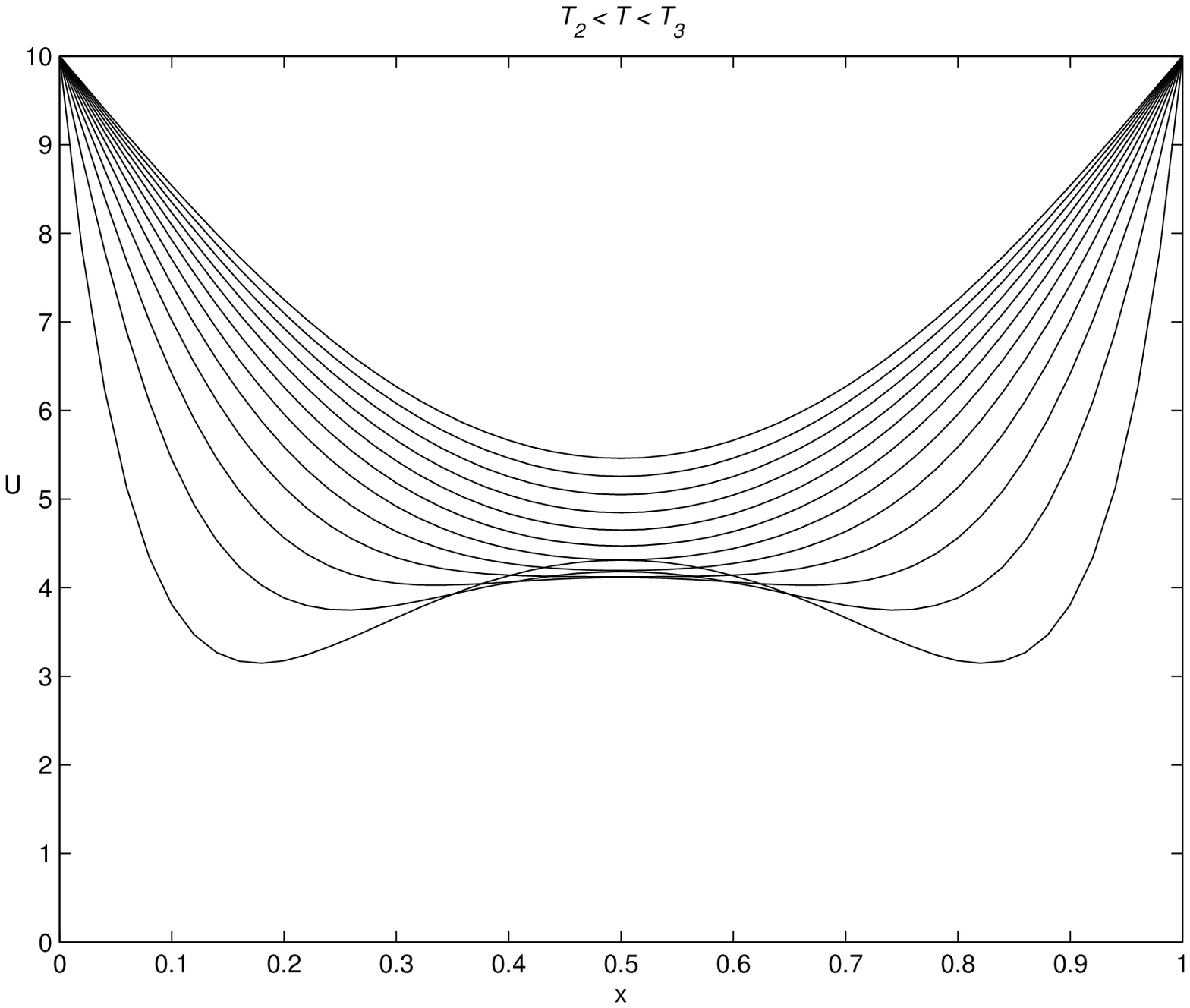}} \caption{}
\end{figure}

\begin{figure}
\scalebox{0.6}{\includegraphics{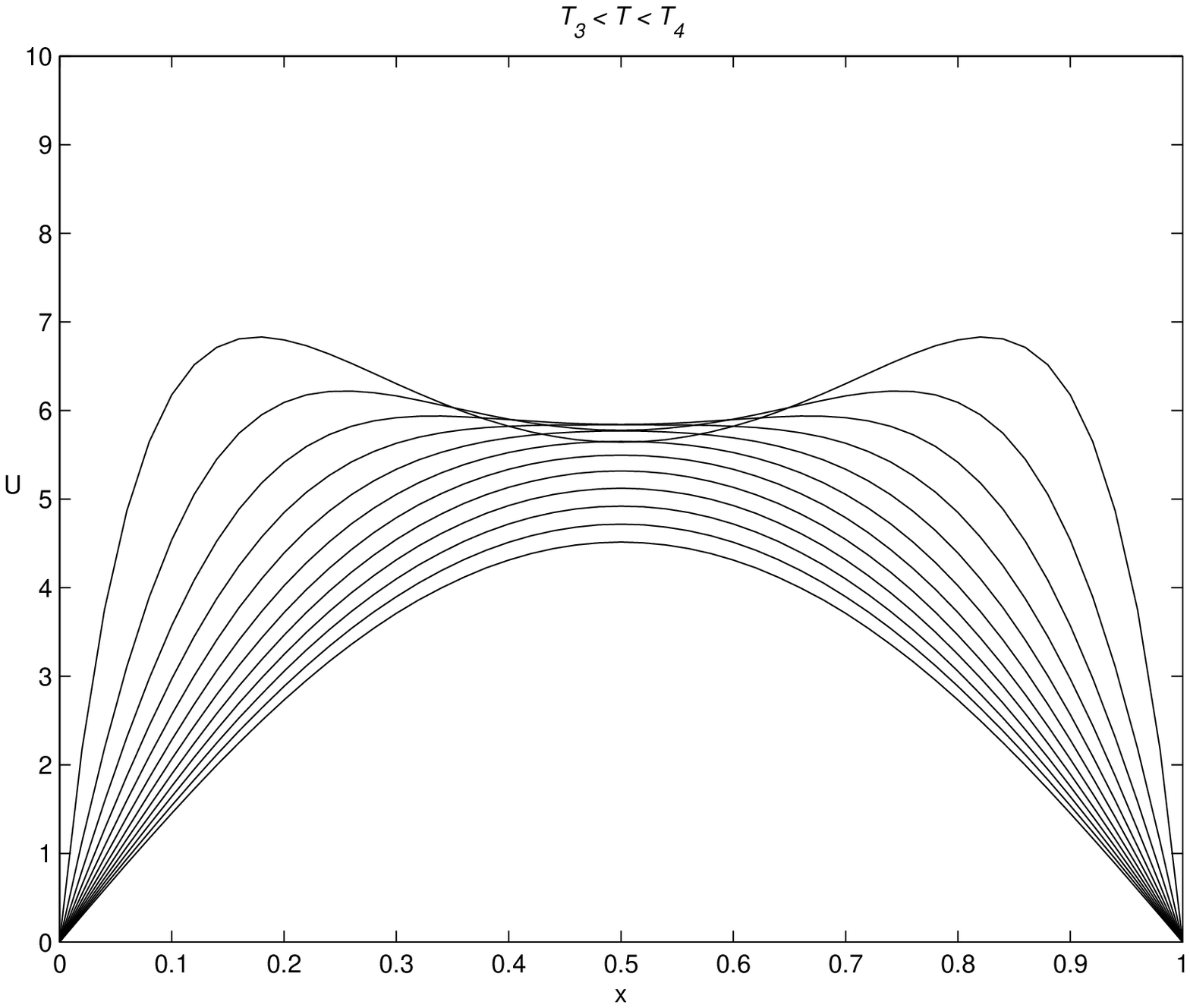}} \caption{}
\end{figure}

\begin{figure}
\scalebox{0.6}{\includegraphics{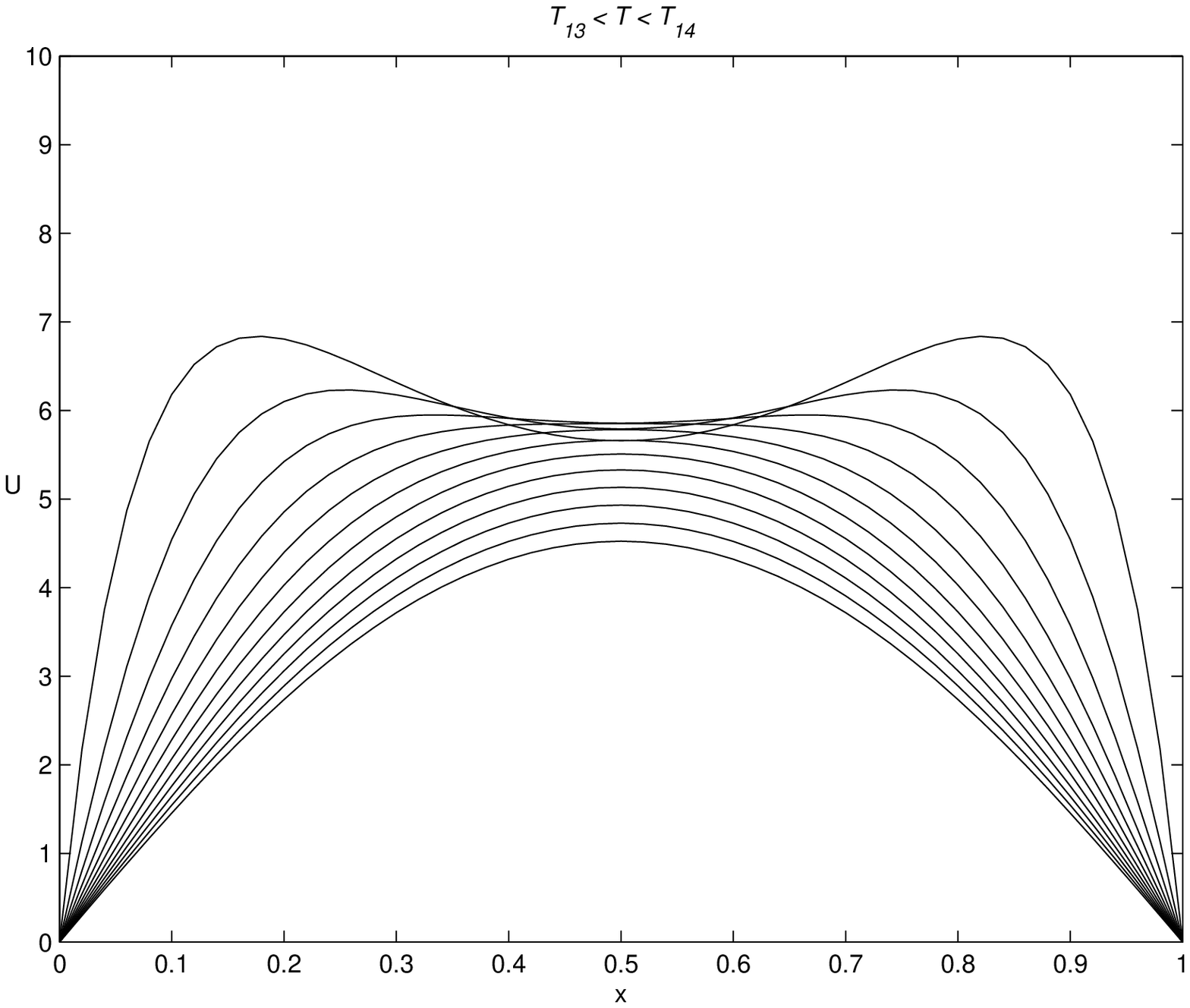}} \caption{}
\end{figure}

\begin{figure}
\scalebox{0.6}{\includegraphics{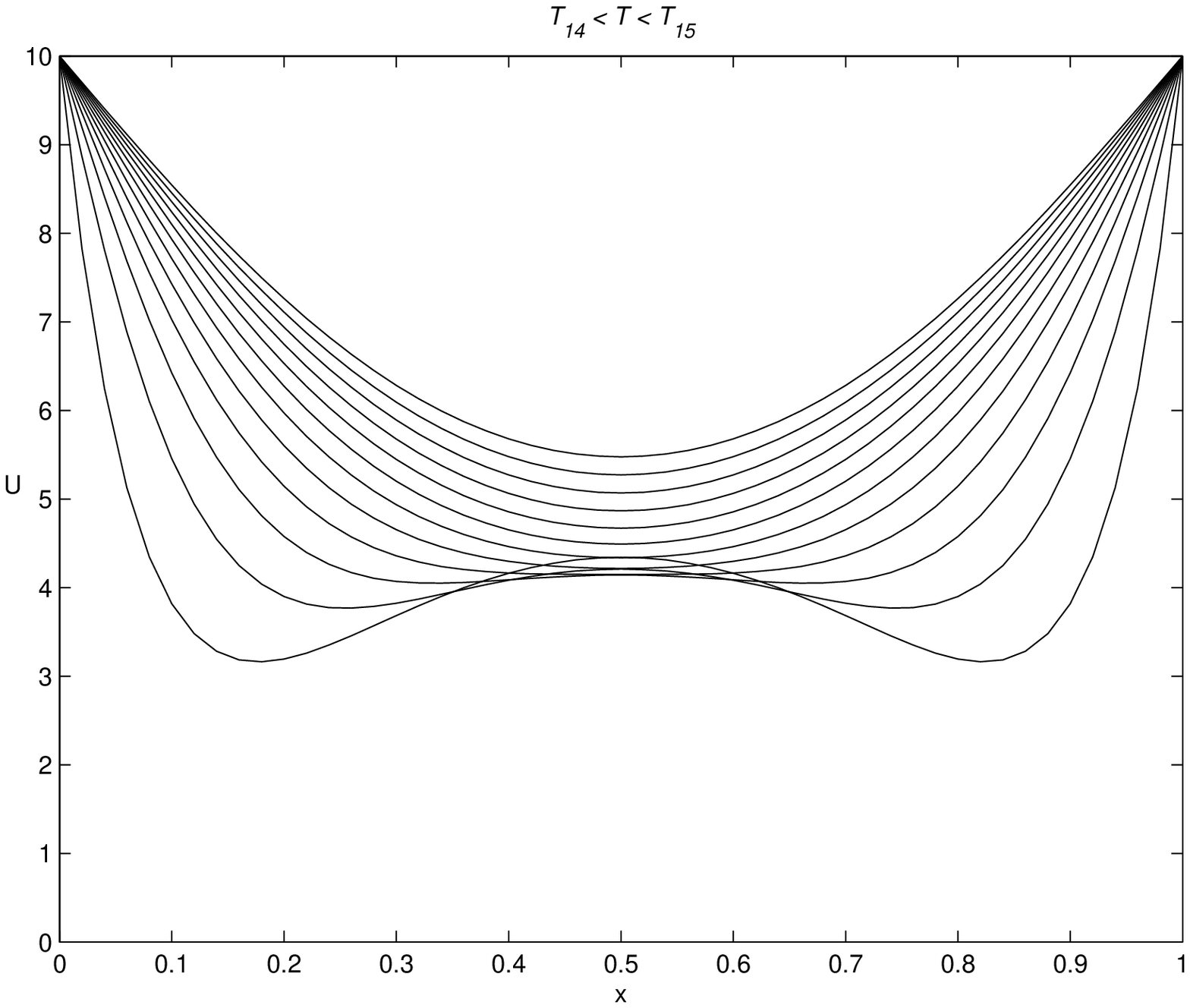}} \caption{}
\end{figure}

\begin{figure}
\scalebox{0.6}{\includegraphics{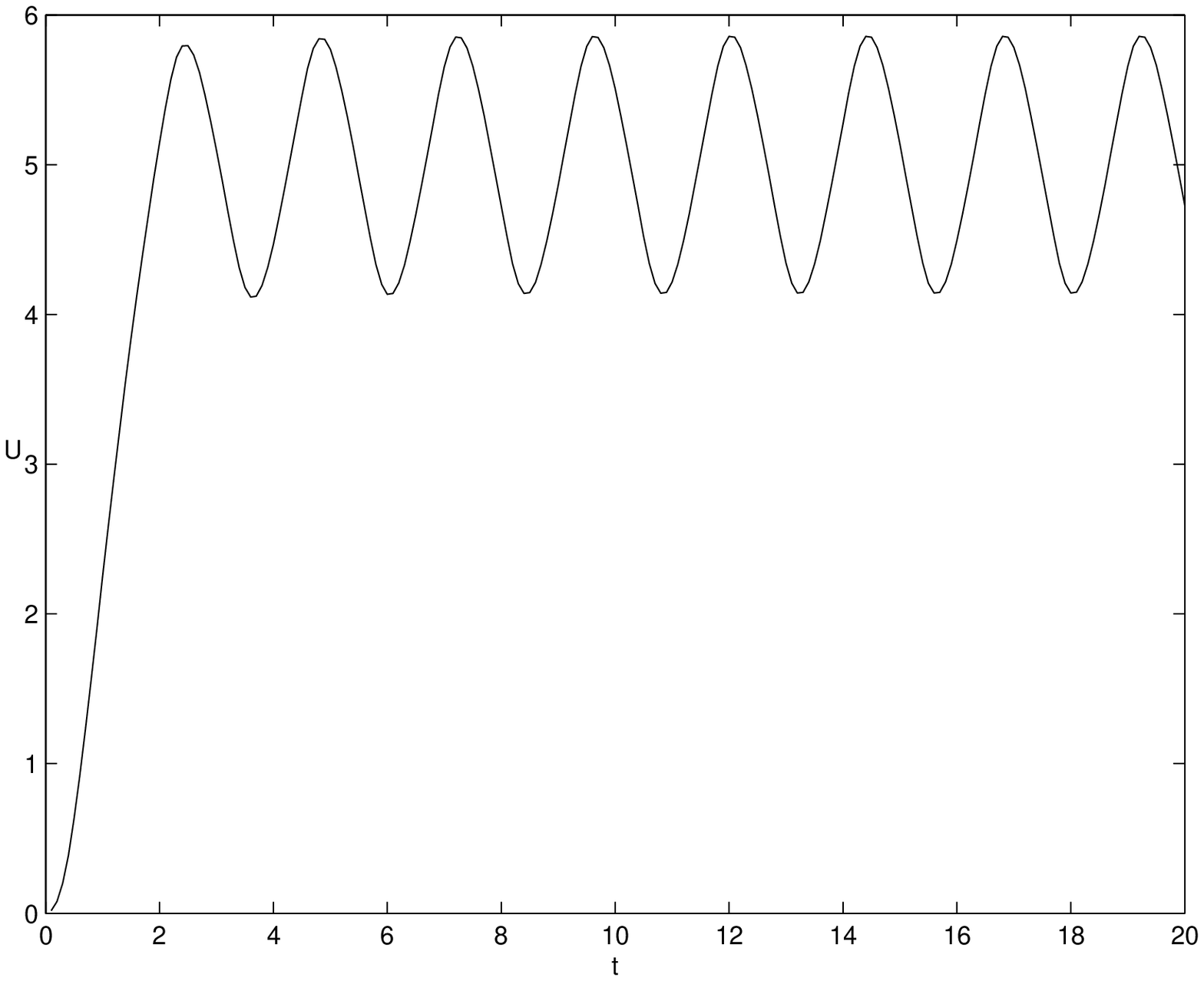}} \caption{}
\end{figure}

\begin{figure}
\scalebox{0.6}{\includegraphics{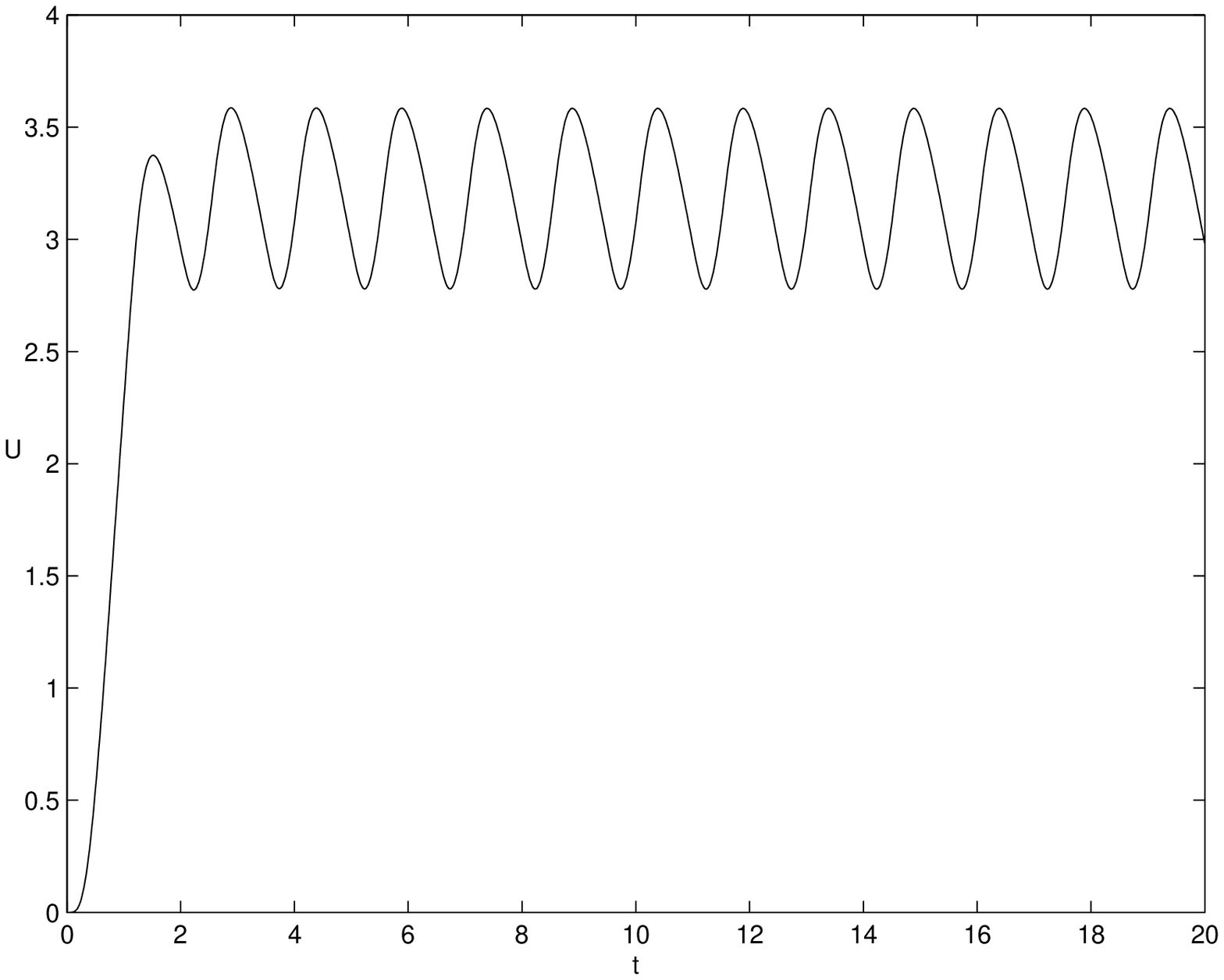}} \caption{}
\end{figure}

\end{document}